%% file: schurQ.tex
\documentclass[11pt]{article}

\usepackage{amsmath}
\usepackage{amsthm}
\usepackage{amssymb}
\usepackage{amscd}

\theoremstyle{definition}
\newtheorem{theorem}{Theorem}[section]
\newtheorem{prop}[theorem]{Proposition}
\newtheorem{lemma}[theorem]{Lemma}
\newtheorem{corollary}[theorem]{Corollary}

\newtheorem{remark}[theorem]{Remark}

\numberwithin{equation}{section}
\newenvironment{demo}[1]{%
  \trivlist
  \item[\hskip\labelsep
        {\it #1.}]
}{%
\hfill\qedsymbol
  \endtrivlist
}

\newcommand\Nat{\mathbb{N}}

\newcommand\Pf{\operatorname{Pf}}
\newcommand\sgn{\operatorname{sgn}}

\newcommand\vectx{\boldsymbol{x}}
\newcommand\vecty{\boldsymbol{y}}
\newcommand\trans{{}^t\!}
\newcommand\vectone{\boldsymbol{1}}
\newcommand\vectzero{\boldsymbol{0}}
\newcommand\ep{\varepsilon}
\newcommand\inv{\operatorname{inv}}

\renewcommand\tilde{\widetilde}

\allowdisplaybreaks

\title{
Pfaffian Formulas and Schur $Q$-Function Identities
}

\author{
Soichi OKADA
\footnote{
Graduate School of Mathematics, Nagoya University, 
Furo-cho, Chikusa-ku, Nagoya 464-8602, Japan, 
{\tt okada@math.nagoya-u.ac.jp}
}
\footnote{
This work was partially supported by 
JSPS Grants-in-Aid for Scientific Research No.~24340003 and No.~15K13425.
}
}

\date{
}

\begin{document}

\maketitle

\begin{abstract}
We establish Pfaffian analogues of the Cauchy--Binet formula 
and the Ishikawa--Wakayama minor-summation formula.
Each of these Pfaffian analogues expresses a sum of products of subpfaffians of two skew-symmetric matrices 
in terms of a single Pfaffian.
By using these Pfaffian formulas we give new transparent proofs to several identities 
for Schur $Q$-functions.
\end{abstract}

\input schurQ1.tex 
\input schurQ2.tex 
\input schurQ3.tex 
\input schurQ4.tex 
\input schurQ5.tex 
\input schurQ6.tex 
\input schurQ7.tex 
\input schurQr.tex 

\end{document}

%% file: schurQ1.tex
\section{%
Introduction
}

The aim of this article is twofold:
Firstly we establish Pfaffian analogues 
of the Cauchy--Binet formula and the Ishikawa--Wakayama minor summation formula \cite{IW} 
for determinants.
Secondly we give new transparent proofs to Schur $Q$-function identities 
by applying general formula for Pfaffians such as these Pfaffian analogues.

Schur $Q$-functions are a family of symmetric functions introduced by Schur \cite{Schur} 
in his study on the projective representations of symmetric groups. 
Schur $Q$-functions play the same role as Schur functions for the linear representation of symmetric groups.
Later Hall (unpublished) and Littlewood \cite{Littlewood} introduce a family of symmetric functions with parameter $t$, 
as a common generalization of Schur functions (the $t=0$ case) and Schur $Q$-functions (the $t=-1$ case).

Schur $Q$-functions appear in various situations parallel to Schur functions:
the projective representations of symmetric groups \cite{Schur},
the cohomology of Lagrangian or orthogonal Grassmannians \cite{Jozefiak2, Pragacz},
the representations of the queer Lie super algebra $\mathfrak{q}(n)$ \cite{Sergeev}, 
the BKP hierarchy \cite{You}.
Also Schur $Q$-functions are expressed as multivariate generating functions 
of shifted tableaux.

In this paper we adopt Nimmo's formula \cite[(A13)]{Nimmo} as a definition of Schur $P$- and $Q$-functions.
This formula is an analogue of the bialternant definition of Schur functions.
Let $\vectx = (x_1, \dots, x_n)$ be a sequence of $n$ indeterminates.
We put
\begin{equation}
\label{eq:matrixA}
A(\vectx) = \left(
 \frac{x_j-x_i}{x_j+x_i}
\right)_{1 \le i, j \le n},
\quad\text{and}\quad
D(\vectx) = \prod_{1 \le i < j \le n}
 \frac{x_j-x_i}{x_j+x_i}.
\end{equation}
For a sequence $\alpha = (\alpha_1, \dots, \alpha_l)$ of nonnegative integers of length $l$, 
let $V_\alpha(\vectx)$ and $W_\alpha(\vectx)$ be the $n \times l$ matrices given by
\begin{equation}
\label{eq:matrixVW}
V_\alpha(\vectx) =
\Bigl(
 x_i^{\alpha_j}
\Bigr)_{1 \le i \le n, 1 \le j \le l},
\quad\text{and}\quad
W_\alpha(\vectx) =
\Bigl(
 \chi(\alpha_j) x_i^{\alpha_j}
\Bigr)_{1 \le i \le n, 1 \le j \le l},
\end{equation}
where $\chi(r) = 2$ if $r > 0$ and $1$ if $r = 0$.
A \emph{strict partition} of length $l$ is 
a strictly decreasing sequence $\lambda = (\lambda_1, \dots, \lambda_l)$ of positive integers.
We write $l = l(\lambda)$.
We define the \emph{Schur $P$-function} $P_\lambda(\vectx)$ and 
the \emph{Schur $Q$-function} $Q_\lambda(\vectx)$ corresponding to a strict partition $\lambda$ 
by putting
\begin{align}
\label{eq:NimmoP}
P_\lambda(\vectx)
 &=
\frac{1}{D(\vectx)} 
\Pf \begin{pmatrix}
 A(\vectx) & V_\alpha(\vectx) \\
 -\trans V_\alpha(\vectx) & O
\end{pmatrix},
\\
\label{eq:NimmoQ}
Q_\lambda(\vectx)
 &=
\frac{1}{D(\vectx)} 
\Pf \begin{pmatrix}
 A(\vectx) & W_\alpha(\vectx) \\
 -\trans W_\alpha(\vectx) & O
\end{pmatrix},
\end{align}
where $\alpha = (\lambda_1, \dots, \lambda_l)$ if $n+l$ is even, 
or $\alpha = (\lambda_1, \dots, \lambda_l, 0)$ if $n+l$ is odd.
Note that $P_\lambda(\vectx) = Q_\lambda(\vectx) = 0$ if $l >n$.

Many of Schur function identities are easily proved by applying determinant formulas.
However some of the known proofs of $Q$-function identities are quite different 
from the proofs of similar Schur function identities.
For example, the Cauchy identity for Schur functions
\begin{equation}
\label{eq:S-Cauchy}
\sum_\lambda s_\lambda(\vectx) s_\lambda(\vecty)
 =
\prod_{i,j} \frac{ 1 }{ 1 - x_i y_j }
\end{equation}
can be proved by using the Cauchy--Binet formula for determinants and the evaluation of Cauchy determinant 
(see \cite[I.4, Example~6]{Macdonald}).
On the other hand, no such direct proof is known for the Cauchy-type identity for Schur $Q$-functions 
\begin{equation}
\label{eq:Q-Cauchy}
\sum_\lambda P_\lambda(\vectx) Q_\lambda(\vecty)
 =
\prod_{i,j} \frac{ 1 + x_i y_j }{ 1 - x_i y_j }.
\end{equation}
See \cite[Abschnitt IV]{Schur}, \cite[\S~4B]{Jozefiak1}, \cite[III.8]{Macdonald} and \cite[Chapter~7]{HH} for 
algebraic proofs.
One of our motivations is to give an elementary linear algebraic proof to $Q$-function identities.

One of the main results of this paper is the following Pfaffian analogue of the Cauchy--Binet formula, 
which reduces to the Cauchy--Binet formula for determinants by specializing $A = O$ and $B = O$.

\begin{theorem}
(Theorem~\ref{thm:Pf-CB} below)
Let $m$ and $n$ be nonnegative integers with the same parity, .
Let $A$ and $B$ be $m \times m$ and $n \times n$ skew-symmetric matrices,
and let $S$ and $T$ be $m \times l$ and $n \times l$ matrices.
Then we have
\begin{multline*}
\sum_K 
 \Pf \begin{pmatrix} A & S([m];K) \\ -\trans S([m];K) & O \end{pmatrix} 
 \Pf \begin{pmatrix} B & T([n];K) \\ -\trans T([n];K) & O \end{pmatrix} 
\\
 =
(-1)^{ \binom{n}{2} }
\Pf \begin{pmatrix}
 A & S \trans T \\
 - T \trans S & -B
\end{pmatrix},
\end{multline*}
where $K$ runs over all subsets of $[l]$ with $\# K \equiv m \equiv n \bmod 2$.
(See Section~2 for notations.)
\end{theorem}

We can give a simple and direct proof to the Cauchy-type identity (\ref{eq:Q-Cauchy}) 
by using this Pfaffian version of the Cauchy--Binet formula as well as the evaluation of Schur Pfaffian.
Also we can use a variant to prove the Pragacz--J\'ozefiak--Nimmo identity for skew $Q$-functions 
\cite{JP, Nimmo}.
In a forthcoming paper, we take this linear algebraic approach to study 
generalizations of Schur $P$- and $Q$-functions such as 
Ivanov's factorial $P$- and $Q$-functions \cite{Ivanov1, Ivanov2} 
and the case $t=-1$ of Hall--Littlewood polynomials associated to the classical root systems \cite{Macdonald00}.

This paper is organized as follows.
After reviewing basic properties of Pfaffians in Section~2, 
we give Pfaffian analogues of the Cauchy--Binet formula 
and the Ishikawa--Wakayama minor-summation formula in Section~3.
In Section~4, we apply the Pfaffian analogue of the Sylvester formula to recover 
Schur's original definition of $Q$-functions from Nimmo's formula.
In Section~5, we give a proof of the Cauchy-type formula for $Q$-functions 
by using the Pfaffian analogue of the Cauchy--Binet formula.
Section~6 is devoted to a linear algebraic proof of the Pragacz--J\'ozefiak--Nimmo formula 
for skew $Q$-functions. 
In Section~7 we use the Pfaffian analogue of the Ishikawa--Wakayama formula 
to derive a Littlewood-type formula for $Q$-functions.

%% file: schurQ2.tex
\section{%
Pfaffians
}

In this section we review basic properties of Pfaffians 
and give a Laplace-type expansion formula.

\subsection{%
Basic properties of Pfaffians
}

Recall the definition and some properties of Pfaffians. (See \cite {IO} for some expositions)
Let $X = \bigl( x_{ij} \bigr)_{1 \le i, j \le 2m}$ be a skew-symmetric matrix of order $2m$.
The \emph{Pfaffian} of $X$, denoted by $\Pf (X)$, is defined by
\begin{equation}
\label{eq:def-Pf}
\Pf (X)
 =
\sum_{\sigma \in F_{2m}} \sgn(\sigma) \prod_{i=1}^m x_{\sigma(2i-1), \sigma(2i)},
\end{equation}
where $F_{2m}$ is the set of permutations $\sigma \in S_{2m}$ satisfying 
$\sigma(1) < \sigma(3) < \dots < \sigma(2m-1)$ and $\sigma(2i-1) < \sigma(2i)$ for $1 \le i \le m$.
Such permutations are in one-to-one correspondence with set-partitions $\pi$ of $\{ 1, 2, \dots, 2m \}$ 
into $m$ disjoint 2-element subsets.
If $\sigma \in F_{2m}$ corresponds to a set-partition 
$\pi = \{ \{ i_1, j_1 \}, \dots, \{ i_m, j_m \} \}$ with $i_k < j_k$ for $1 \le k \le m$, 
then we have
$$
\sgn(\sigma) = (-1)^{\inv(i_1, j_1, \dots, i_m, j_m)},
$$
where $\inv(\alpha_1, \dots, \alpha_{2m})$ is the number of pairs $(k,l)$ such that 
$k < l$ and $\alpha_k > \alpha_l$.
Note that the right hand side is independent of the ordering of blocks of $\pi$.
Since $m \equiv \binom{2m}{2} \bmod 2$, it follows from the definition of Pfaffians (\ref{eq:def-Pf}) that
\begin{equation}
\label{eq:Pf(-X)}
\Pf (-X) = (-1)^{\binom{2m}{2}} \Pf X.
\end{equation}

Pfaffians are multilinear in the following sense.
Let $X = \bigl( x_{ij} \bigr)_{1 \le i, j \le n}$ be a skew-symmetric matrix 
and fix a row/column index $k$.
If the entries of the $k$th row and $k$th column of $X$ are written as 
$x_{i,j} = \alpha x'_{i,j} + \beta x''_{i,j}$ for $i=k$ or $j=k$, 
then
$$
\Pf X = \alpha \Pf X' + \beta \Pf X'',
$$
where $X'$ (resp. $X''$) is the skew-symmetric matrix obtained from $X$ by 
replacing the entries $x_{ij}$ for $i=k$ or $j=k$ with $x'_{ij}$ (resp. $x''_{ij}$).

If $X$ is an $n \times n$ skew-symmetric matrix and $U$ is an $n \times n$ matrix, then we have
\begin{equation}
\label{eq:UXU}
\Pf \bigl( \trans U X U \bigr) = \det (U) \Pf (X).
\end{equation}
It follows that Pfaffians are alternating, i.e., 
if $\sigma \in S_n$, we have
$$
\Pf \bigl( x_{\sigma(i),\sigma(j)} \bigr)_{1 \le i, j \le n}
 =
\sgn(\sigma) \Pf \bigl( x_{i,j} \bigr)_{1 \le i, j \le n}.
$$
Also we see that, if $Y$ is the skew-symmetric matrix obtained from $X$ 
by adding the $k$th row multiplied by a scalar $\alpha$ to the $l$th row and then 
adding the $k$th column multiplied by $\alpha$ to the $l$th column,
the we have $\Pf Y = \Pf X$.

We use the following notations for submatrices.
For a positive integer $n$, we put $[n] = \{ 1, 2, \dots, n \}$.
Given a subset $I \subset [n]$, we put $\Sigma(I) = \sum_{i \in I} i$.
For an $M \times N$ matrix $X = \bigl( x_{i,j} \bigr)_{1 \le i \le M, 1 \le j \le N}$ 
and subsets $I \subset [M]$ and $J \subset [N]$, 
we denote by $X(I;J)$ the submatrix of $X$ obtained by picking up rows indexed by $I$ and 
columns indexed by $J$.
If $X$ is a skew-symmetric matrix, then we write $X(I)$ for $X(I;I)$.
We use the convention that $\det X(\emptyset;\emptyset) = 1$ and $\Pf X(\emptyset) = 1$.

For an $n \times n$ skew-symmetric matrix $X = \bigl( x_{i,j} \bigr)_{1 \le i, j \le n}$, 
we have the following expansion formula along the $k$th row/column:
\begin{equation}
\label{eq:Pf-expansion}
\Pf X
 = 
\sum_{i=1}^{k-1} (-1)^{k+i-1} x_{i,k} \Pf X([n] \setminus \{ i, k \})
 +
\sum_{i=k+1}^n (-1)^{k+i-1} x_{k,i} \Pf X([n] \setminus \{ k, i \}).
\end{equation}

Knuth \cite{Knuth} gave the following Pfaffian analogue of the Sylvester identity for determinant.

\begin{prop}
\label{prop:Pf-Sylvester}
(Knuth \cite[(2.5)]{Knuth})
Let $n$ and $l$ be even integers 
and let $X$ be an $(n+l) \times (n+l)$ skew-symmetric matrix.
\begin{enumerate}
\item[(1)]
If $\Pf X([n]) \neq 0$, then we have
\begin{equation}
\label{eq:Pf-Sylvester1}
\Pf \left(
 \frac{ \Pf X([n] \cup \{ n+i, n+j \}) }
      { \Pf X([n]) }
\right)_{1 \le i, j \le l}
 =
\frac{ \Pf X }
     { \Pf X([n]) }.
\end{equation}
\item[(2)]
If $\Pf X([l+1,l+n]) \neq 0$, then we have
\begin{equation}
\label{eq:Pf-Sylvester2}
\Pf \left(
 \frac{ \Pf X(\{ i, j \} \cup [l+1,l+n]) }
      { \Pf X([l+1,l+n]) }
\right)_{1 \le i, j \le l}
 =
\frac{ \Pf X }
     { \Pf X([l+1,l+n]) },
\end{equation}
where $[l+1,l+n] = \{ l+1, l+2, \dots, l+n \}$.
\end{enumerate}
\end{prop}

The following evaluation formula of Schur Pfaffian is useful in various places of this paper.

\begin{prop}
\label{prop:Pf-Schur}
(Schur \cite[p.~226]{Schur}, see also \cite[Section~4]{Knuth})
Let $n$ be an even integer.
For a sequence $\vectx = (x_1, \dots, x_n)$ of variables, we have
\begin{equation}
\label{eq:Pf-Schur}
\Pf \left(
 \frac{x_j-x_i}{x_j+x_i}
\right)_{1 \le i, j \le n}
 =
\prod_{1 \le i < j \le n} \frac{x_j-x_i}{x_j+x_i}.
\end{equation}
\end{prop}

\subsection{%
Laplace-type expansion formulas for Pfaffian
}

The following expansion formula is stated in \cite[(12)]{C} without proof.

\begin{prop}
\label{prop:Pf-Laplace}
Let $m$ and $n$ be nonnegative integers such that $m+n$ is even.
For an $m \times m$ skew-symmetric matrix $Z = \bigl( z_{i,j} \bigr)_{1 \le i, j \le m}$, 
an $n \times n$ skew-symmetric matrix $Z' = \bigl( z'_{i,j} \bigr)_{1 \le i, j \le n}$, 
and an $m \times n$ matrix $W = \bigl( w_{i,j} \bigr)_{1 \le i \le m, 1 \le j \le n}$, 
we have
\begin{equation}
\label{eq:Pf-Laplace}
\Pf \begin{pmatrix}
 Z & W \\
 -\trans W & Z'
\end{pmatrix}
 =
\sum_{I,J}
 \ep(I,J)
 \Pf Z(I) \Pf Z'(J) \det W( [m] \setminus I ; [n] \setminus J ),
\end{equation}
where the sum is taken over all pairs of even-element subsets $(I,J)$ such that 
$I \subset [m]$, $J \subset [n]$ and $m - \# I = n - \# J$, 
and the coefficient $\ep(I,J)$ is given by
$$
\ep(I,J)
 = (-1)^{ \Sigma(I) + \Sigma(J) + \binom{m}{2} + \binom{n}{2} + \binom{k}{2} },
\quad
k = m - \# I = n - \# J.
$$
\end{prop}

If $m=1$, then the formula (\ref{eq:Pf-Laplace}) reduces to the expansion formula (\ref{eq:Pf-expansion}) 
along the first row/column.

\begin{demo}{Proof}
We put $[n]' = \{ 1' , 2', \dots, n' \}$ and 
label the rows and columns of $\begin{pmatrix} Z & W \\ - \trans W & Z' \end{pmatrix}$ 
by $[m] \sqcup [n]' = \{ 1, 2, \dots, m, 1', 2', \dots, n' \}$ 
with $1 < 2 < \dots < m < 1' < 2' < \dots < n'$.
For an even-element subset $I$ of $[m] \sqcup [n]'$, we denote by $F_I$ 
the set of all set-partitions of $I$ into $2$-element subsets.
Given a partition $\pi \in F_{[m] \sqcup [n]'}$, we put
$$
\pi_i = \{ b \in \pi : \# (b \cap [m]) = i \}
\quad\text{for $i = 0, 1, 2$.}
$$
Then there are subsets $I \subset [m]$ and $J' \subset [n]'$ such that
$\pi_2 \in F_I$, $\pi_0 \in F_{J'}$ and $m - \# I = n - \# J'$.
Moreover, if $[m] \setminus I = \{ r_1, \dots, r_k \}$ and 
$[n]' \setminus J' = \{ s'_1, \dots, s'_k \}$ with $r_1 < \dots < r_k$, 
$s_1 < \dots < s_k$, then there exists a unique permutation $\sigma \in S_k$ 
such that $\pi_1 = \bigl\{ \{ r_1, s'_{\sigma(1)} \}, \dots, \{ r_k, s'_{\sigma'(k)} \} \bigr\}$.
The correspondence $\pi \mapsto (\pi_2, \sigma, \pi_0)$ gives a bijection 
$$
F_{[m] \sqcup [n]'} \to \bigsqcup_{(I,J)} F_I \times S_k \times F_{J'},
$$
where $(I,J)$ runs over all pairs of even-element subsets $I \subset [m]$ 
and $J \subset [n]$ such that $m - \# I = n - \# J$, 
and $J' = \{ j' : j \in J \}$.
Let $\pi_2 = \bigl\{ \{ p_1, p_2 \}, \dots, \{ p_{2(m-k)-1}, p_{2(m-k)} \} \bigr\}$ 
and $\pi_0 = \bigl\{ \{ q'_1, q'_2 \}, \dots, \{ q'_{2(n-k)-1}, q'_{2(n-k)} \} \bigr\}$ 
with $p_{2i-1} < p_{2i}$ and $q_{2j-1} < q_{2j}$.
Then the inversion number of the permutation associated to $\pi$ is given by 
\begin{align*}
&
\inv \bigl( 
 p_1, p_2, \dots, p_{2(m-k)},
 r_1, s'_{\sigma(1)}, \dots, r_k, s'_{\sigma(k)}, 
 q'_1, q'_2, \dots, q'_{2(n-k)}
\bigr)
\\
 &
=
\inv \bigl( p_1, p_2, \dots, p_{2(m-k)} \bigr)
+
\inv \bigl( r_1, s'_{\sigma(1)}, \dots, r_k, s'_{\sigma(k)} \bigr)
+
\inv \bigl( q'_1, q'_2, \dots, q'_{2(n-k)} \bigr)
\\
 &
 \quad
+
\# \{ (i,j) \in I \times ( [m] \setminus I ) : i > j \}
 +
\binom{k}{2}
 +
\# \{ (i',j') \in ([n]' \setminus J') \times J : i' > j' \}.
\end{align*}
Since $r_1 < \dots < r_k < s'_1 < \dots < s'_k$, we have
$$
\inv \bigl( r_1, s'_{\sigma(1)}, \dots, r_k, s'_{\sigma(k)} \bigr)
 =
\inv \bigl( \sigma(1), \dots, \sigma(k) \bigr).
$$
Also we have
\begin{gather*}
\# \{ (i,j) \in I \times ( [m] \setminus I ) : i > j \}
 =
\Sigma(I) - \binom{m-k+1}{2},
\\
\# \{ (i',j') \in ([n]' \setminus J') \times J' : i' > j' \}
 =
k(n-k) + \binom{n-k+1}{2} - \Sigma(J).
\end{gather*}
Since $m \equiv n \equiv k \bmod 2$ by the assumption, we see that
\begin{align*}
\binom{m-k+1}{2} + k(n-k) + \binom{n-k+1}{2}
 &=
\binom{m}{2} + \binom{n}{2} - (k+1)m + m + n
\\
 &\equiv
\binom{m}{2} + \binom{n}{2} \bmod 2.
\end{align*}
Hence we have
\begin{align*}
&
\inv \bigl( 
 p_1, p_2, \dots, p_{2(m-k)},
 r_1, s'_{\sigma(1)}, \dots, r_k, s'_{\sigma(k)}, 
 q'_1, q'_2, \dots, q'_{2(n-k)}
\bigr)
\\
 &
\equiv
\inv \bigl( p_1, p_2, \dots, p_{2(m-k)} \bigr)
+
\inv \bigl( \sigma(1), \dots, \sigma(k) \bigr)
+
\inv \bigl( q'_1, q'_2, \dots, q'_{2(n-k)} \bigr)
\\
 &
 \quad
+ \Sigma(I) + \Sigma(J) + \binom{m}{2} + \binom{n}{2} + \binom{k}{2}.
\end{align*}
Now (\ref{eq:Pf-Laplace}) follows from the definition of Pfaffians (\ref{eq:def-Pf}).
\end{demo}

By considering the case where $Z'$ or $W$ is the zero matrix 
in Proposition~\ref{prop:Pf-Laplace}, we obtain the following corollary.
We denote the $p \times q$ zero matrix by $O_{p,q}$ and 
write simply $O$ for $O_{p,q}$ if there is no confusion on the size.

\begin{corollary}
\label{cor:Pf-Laplace}
Suppose that $m+n$ is even.
\begin{enumerate}
\item[(1)]
If $Z$ is an $m \times m$ skew-symmetric matrix and $W$ is an $m \times n$ matrix, then we have
\begin{equation}
\label{eq:Pf-Laplace1}
\Pf \begin{pmatrix}
 Z & W \\
 -\trans W & O_{n,n}
\end{pmatrix}
 =
\begin{cases}
\displaystyle\sum_{I}
 (-1)^{\Sigma(I) + \binom{m}{2}}
 \Pf Z(I) \det W( [m] \setminus I ; [n] )
 &\text{if $m > n$,}
\\
(-1)^{\binom{m}{2}} \det W
 &\text{if $m = n$,}
\\
0
 &\text{if $m < n$,}
\end{cases}
\end{equation}
where $I$ runs over all $(m-n)$-element subsets of $[n]$.
\item[(2)]
If $Z$ and $Z'$ are $m \times m$ and $n \times n$ skew-symmetric matrices respectively, then we have
\begin{equation}
\label{eq:Pf-Laplace2}
\Pf \begin{pmatrix}
 Z & O_{m,n} \\
 O_{n,m} & Z'
\end{pmatrix}
 =
\begin{cases}
 \Pf Z \cdot \Pf Z' &\text{if $m$ and $n$ are even,} \\
 0 &\text{otherwise}.
\end{cases}
\end{equation}
\end{enumerate}
\end{corollary}

\begin{demo}{Proof}
(1)
If $Z' = O$, then we have $\Pf Z'(J) = 0$ unless $J = \emptyset$.

(2)
If $W = O$, then we have $\det W([m] \setminus I; [n] \setminus J) = 0$ unless $I = [m]$ and $J = [n]$.
\end{demo}

%% file: schurQ3.tex
\section{%
Cauchy--Binet type Pfaffian formulas
}

In this section we give Pfaffian analogues of the Cauchy--Binet formula 
and the Ishikawa--Wakayama minor-summation formula \cite{IW}.
These are our main results of this paper.

First we consider the following special case of Proposition~\ref{prop:Pf-Laplace}.

\begin{lemma}
\label{lem:Pf-CBIW}
Let $m$, $n$ and $l$ be nonnegative integers with $m \equiv n \bmod 2$.
We put
$$
E^{(m+l,n+l)}_l
 =
\begin{pmatrix}
 O_{m,n} & O_{m,l} \\
 O_{l,n} & E_l
\end{pmatrix},
$$
where $E_l$ is the $l \times l$ identity matrix.
If $Z$ and $Z'$ are $(m+l) \times (m+l)$ and $(n+l) \times (n+l)$ skew-symmetric matrices respectively, 
then we have
\begin{multline}
\label{eq:Pf-CBIW}
\Pf \begin{pmatrix}
 Z & E^{(m+l,n+l)}_l \\
 -\trans E^{(m+l,n+l)}_l & Z'
\end{pmatrix}
\\
 =
\sum_K
 (-1)^{ \binom{l - \# K}{2} }
 \Pf Z([m] \cup (m+K)) \Pf Z'([n] \cup (n+K)),
\end{multline}
where $K$ runs over all subsets of $[l]$ with $\# K \equiv m \bmod 2$ and
 $m+K = \{ m+k : k \in K \}$, $n+K = \{ n+k : k \in K \}$.
\end{lemma}

\begin{demo}{Proof}
We substitute $W = E^{(m+l,n+l)}_l$ in Proposition~\ref{prop:Pf-Laplace}.
Let $I$ and $J$ be even-element subsets of $[m+l]$ and $[n+l]$ respectively 
such that $m+ l - \# I = n + l - \# J$.
If $[m] \not\subset I$ or $[n] \not\subset J$, then we have 
$\det W([m+l] \setminus I; [n+l] \setminus J) = 0$.
If $[m] \subset I$ and $[n] \subset J$, then we can write 
$I = [m] \sqcup (m+I')$ and $J = [n] \sqcup (n+J')$ for some subsets $I'$, $J' \subset [l]$,
and we see that
$$
\det W([m+l] \setminus I;[n+l] \setminus J)
 = \det E_l ([l] \setminus I';[l] \setminus J')
 =
\begin{cases}
 1 &\text{if $I' = J'$,} \\
 0 &\text{otherwise.}
\end{cases}
$$
Hence $\det W([m+l] \setminus I;[n+l] \setminus J) = 0$ 
unless $I = [m] + (m+K)$ and $J = [n] + (n+K)$ for some subset $K \subset [l]$.
In this case,
$$
\Sigma(I) = \binom{m}{2} + m + m \# K + \Sigma(K),
\quad
\Sigma(J) = \binom{n}{2} + n + n \# K + \Sigma(K),
$$
and
$$
\Sigma(I) + \Sigma(J) + \binom{m+l}{2} + \binom{n+l}{2} + \binom{(m+l)-(m+\# K)}{2}
 \equiv
\binom{l-\# K}{2} \bmod 2.
$$
This completes the proof.
\end{demo}

We use Lemma~\ref{lem:Pf-CBIW} to derive a Pfaffian analogue of the Cauchy--Binet formula.

\begin{theorem}
\label{thm:Pf-CB}
Let $m$ and $n$ be nonnegative integers with the same parity, .
Let $A$ and $B$ be $m \times m$ and $n \times n$ skew-symmetric matrices,
and let $S$ and $T$ be $m \times l$ and $n \times l$ matrices.
Then we have
\begin{multline}
\label{eq:Pf-CB1}
\sum_K (-1)^{ \binom{\# K}{2} }
 \Pf \begin{pmatrix} A & S([m];K) \\ -\trans S([m];K) & O \end{pmatrix} 
 \Pf \begin{pmatrix} B & T([n];K) \\ -\trans T([n];K) & O \end{pmatrix} 
\\
 =
\Pf \begin{pmatrix}
 A & S \trans T \\
 - T \trans S & B
\end{pmatrix},
\end{multline}
\begin{multline}
\label{eq:Pf-CB2}
\sum_K 
 \Pf \begin{pmatrix} A & S([m];K) \\ -\trans S([m];K) & O \end{pmatrix} 
 \Pf \begin{pmatrix} B & T([n];K) \\ -\trans T([n];K) & O \end{pmatrix} 
\\
 =
(-1)^{ \binom{n}{2} }
\Pf \begin{pmatrix}
 A & S \trans T \\
 - T \trans S & -B
\end{pmatrix},
\end{multline}
where $K$ runs over all subsets of $[l]$ with $\# K \equiv m \equiv n \bmod 2$.
\end{theorem}

\begin{remark}
It follows from (\ref{eq:Pf-Laplace1}) that both formulas (\ref{eq:Pf-CB1}) and (\ref{eq:Pf-CB2}) 
reduce to the Cauchy--Binet formula for determinants if we put $A = 0$ and $B = 0$:
\begin{equation}
\label{eq:CB}
\sum_K \det S([m];K) \det T([m];K)
 =
\det (S \trans T),
\end{equation}
where $K$ runs over all $m$-element subsets.
\end{remark}

\begin{demo}{Proof}
Apply Lemma~\ref{lem:Pf-CBIW} to the matrices
$$
Z = \begin{pmatrix}
 A & S \\
 -\trans S & O
\end{pmatrix}
\quad\text{and}\quad
Z' = \begin{pmatrix}
 B & -T \\
 \trans T & O
\end{pmatrix}.
$$
Then we have
\begin{gather*}
\Pf Z([m]+(m+K))
 =
\Pf \begin{pmatrix}
 A & S([m];K) \\
 -\trans S([m];K) & O
\end{pmatrix},
\\
\Pf Z'([n]+(n+K))
 =
(-1)^{\# K}
\Pf \begin{pmatrix}
 B & T([n];K) \\
 - \trans T([n];K) & O
\end{pmatrix}.
\end{gather*}
We compute the Pfaffian on the right hand side of (\ref{eq:Pf-CBIW}).
By using the relation (\ref{eq:UXU}) with
$$
X = \begin{pmatrix}
 A & S & O & O \\
 -\trans S & O & O & E_l \\
 O & O & B & -T \\
 O & -E_l & \trans T & O
\end{pmatrix}
\quad\text{and}\quad
U = \begin{pmatrix}
 E_m & O & O & O \\
 O & \trans T & E_l & O \\
 O & E_n & O & O \\
 \trans S & O & O & E_l
\end{pmatrix},
$$
and then by using Corollary~\ref{cor:Pf-Laplace}, we see that
$$
(-1)^{nl} \Pf X
 =
\Pf \begin{pmatrix}
 A & S \trans T & O & O \\
 - T \trans S & B & O & O \\
 O & O & O & E \\
 O & O & -E & O
\end{pmatrix}
 =
(-1)^{\binom{l}{2}}
\Pf \begin{pmatrix}
 A & S \trans T \\
 - T \trans S & B
\end{pmatrix}.
$$
Therefore we have
\begin{multline*}
\Pf \begin{pmatrix}
 A & S \trans T \\
 - T \trans S & B
\end{pmatrix}
\\
 =
\sum_K
 (-1)^{ \binom{l-\# K}{2} + \# K + nl - \binom{l}{2} }
 \Pf \begin{pmatrix}
  A & S([m];K) \\
  -\trans S([m];K) & O
 \end{pmatrix}
 \Pf \begin{pmatrix}
  B & T([n];K) \\
  -\trans T([n];K) & O
 \end{pmatrix}.
\end{multline*}
Since $\# K \equiv n \bmod 2$, we have
$$
\binom{l-\# K}{2} + \# K + nl - \binom{l}{2}
 =
\binom{\# K}{2} + l(n-\# K) + 2 \# K
 \equiv
\binom{\# K}{2} \bmod 2,
$$
and obtain the desired formula (\ref{eq:Pf-CB1}).

Equation (\ref{eq:Pf-CB2}) is obtained by replacing $B$ with $-B$ in (\ref{eq:Pf-CB1}).
In fact, by multiplying the last $k$ rows/columns by $-1$ 
and then by using (\ref{eq:Pf(-X)}), we have
$$
\Pf \begin{pmatrix}
 -B & T([n];K) \\
 -\trans T([n];K) & O
\end{pmatrix}
 =
(-1)^{ \# K + \binom{n+\# K}{2} }
\Pf \begin{pmatrix}
 B & T([n];K) \\
 - \trans T([n];K) & O
\end{pmatrix}.
$$
Since $\# K \equiv n \bmod 2$, we have
$$
\# K + \binom{n+\# K}{2}
 =
\binom{n}{2} + \binom{\# K}{2} + (n+1) \# K
 \equiv
\binom{n}{2} + \binom{\# K}{2} \bmod 2,
$$
and obtain (\ref{eq:Pf-CB2}).
\end{demo}

Another application of Lemma~\ref{lem:Pf-CBIW} is the following Pfaffian analogue of 
the Ishikawa--Wakayama minor-summation formula.

\begin{theorem}
\label{thm:Pf-IW}
Let $m$ be an even integer and $l$ be a positive integer.
For an $m \times m$ skew-symmetric matrix $A$, 
an $l \times l$ skew-symmetric matrix $B$, 
and an $m \times l$ matrix $S$, we have
\begin{equation}
\label{eq:Pf-IW}
\sum_K 
 \Pf B(K) 
 \Pf \begin{pmatrix}
  A & S([m];K) \\
  - \trans S([m];K) & O
 \end{pmatrix}
=
\Pf \left(
 A - S B \trans S
\right),
\end{equation}
where $K$ runs over all even-element subsets of $[l]$.
\end{theorem}

\begin{remark}
It follows from (\ref{eq:Pf-Laplace1}) that (\ref{eq:Pf-IW}) reduces 
to the minor-summation formula (\cite[Theorem~1]{IW}) if $A = O$:
\begin{equation}
\label{eq:IW}
\sum_K \Pf B(K) \det S([m];K)
 =
\Pf \left( S B \trans S \right),
\end{equation}
where $K$ runs over all $m$-element subsets of $[l]$.
\end{remark}

\begin{demo}{Proof}
We apply Lemma~\ref{lem:Pf-CBIW} (with $n=0$) to the matrices
$$
Z = \begin{pmatrix}
 A & S \\
 -\trans S & O
\end{pmatrix},
\quad
Z' = - B.
$$
Since $\Pf Z'(K) = (-1)^{\binom{\# K}{2}} \Pf B(K)$ by (\ref{eq:Pf(-X)}), we have
$$
\Pf \begin{pmatrix}
 A & S & O \\
 -\trans S & O & E_l \\
 0 & -E_l & -B
\end{pmatrix}
 =
\sum_{K \subset [l]}
 (-1)^{ \binom{l- \# K}{2} + \binom{\# K}{2} }
 \Pf \begin{pmatrix}
  A & S([m];K) \\
  -\trans S([m];K) & O
 \end{pmatrix}
 \Pf B(K).
$$
By using (\ref{eq:UXU}) with
$$
X
 = 
\begin{pmatrix}
 A & S & O \\
 -\trans S & O & E \\
 O & -E & -B
\end{pmatrix},
\quad
U
 =
\begin{pmatrix}
E & O & O \\
- B \trans S & E & O \\
\trans S & O & E
\end{pmatrix}
$$
and then by using Corollary~\ref{cor:Pf-Laplace}, we obtain
$$
\Pf X
 =
\Pf \begin{pmatrix}
A - S B \trans S & O & O \\
O & O & E \\
O & -E & -B
\end{pmatrix}
 =
(-1)^{\binom{l}{2}}
\Pf \left( A - S B \trans S \right).
$$
Hence the proof is completed by using the congruence 
$\binom{l-\# K}{2} + \binom{\# K}{2} + \binom{l}{2} \equiv 0 \bmod 2$.
\end{demo}

\begin{remark}
From Lemma~\ref{lem:Pf-CBIW}, we can derive the following summation formula for Pfaffians 
\cite[Theorem~3]{IW}:
\begin{equation}
\label{eq:IW2}
\sum_{I,J} (-1)^{\binom{l-\# I}{2}} \det T(I;J) \Pf A(I) \Pf B(J) 
 =
\Pf \begin{pmatrix}
 A & E_l \\
 -E_l & T B \trans T
\end{pmatrix},
\end{equation}
where $A$ and $B$ are $l \times l$ skew-symmetric matrices, $T$ is an $l \times l$ matrix, 
and the summation is taken over all pairs of even-element subsets $I$, $J \subset [l]$ 
such that $\# I = \# J$.
In fact, if we consider the case $m=n=0$ of Lemma~\ref{lem:Pf-CBIW}, we obtain
$$
\sum_K (-1)^{\binom{l-\# K}{2}} \Pf Z(K) \Pf Z'(K)
 =
\Pf \begin{pmatrix}
 Z & E_l \\
 -E_l & Z'
\end{pmatrix}.
$$
By taking $Z = A$ and $Z' = T B \trans T$ and using the minor-summation formula (\ref{eq:IW}), 
we obtain (\ref{eq:IW2}).
\end{remark}

%% file: schurQ4.tex
\section{%
Schur's original definition of $Q$-functions
}

In this section, we recover Schur's original definition \cite{Schur} 
of $Q$-functions from Nimmo's formula (\ref{eq:NimmoQ}) 
by applying the Pfaffian analogue of the Sylvester formula (Proposition~\ref{prop:Pf-Sylvester}).
Macdonald \cite[III. 8]{Macdonald} proves Part (3) of the following theorem 
by considering the generating function of Hall--Littlewood functions, 
And Stembridge's derivation \cite[Theorem~6.1]{Stembridge} is based on the combinatorial definition of $Q$-functions 
and the lattice path method.

\begin{theorem}
\label{thm:Schur}
(Schur \cite{Schur})
\begin{enumerate}
\item[(1)]
The generating function of Schur $Q$-functions corresponding to partitions of length $\le 1$ 
is given by
\begin{equation}
\label{eq:Schur1}
\sum_{r \ge 0} Q_{(r)}(\vectx) z^r 
=
\prod_{i=1}^n \frac{ 1 + x_i z }{ 1 - x_i z},
\end{equation}
where $Q_{(0)}(\vectx) = 1$.
\item[(2)]
The generating function of Schur $Q$-functions corresponding to partitions of length $\le 2$ 
is given by
\begin{equation}
\label{eq:Schur2}
\sum_{r, s \ge 0} Q_{(r,s)}(\vectx) z^r w^s
=
\frac{z-w}{z+w}
\left(
 \prod_{i=1}^n \frac{ 1 + x_i z }{ 1 - x_i z}
 \prod_{i=1}^n \frac{ 1 + x_i w }{ 1 - x_i w}
-
1
\right),
\end{equation}
where $Q_{(0,0)}(\vectx) = 0$ and
$$
Q_{(r,s)}(\vectx) = - Q_{(r,s)}(\vectx),
\quad
Q_{(r,0)}(\vectx) = - Q_{(0,r)}(\vectx) = Q_{(r)}(\vectx)
$$
for positive integers $r$ and $s$.
\item[(3)]
For a sequence of nonnegative integers $\alpha = (\alpha_1, \dots, \alpha_l)$, we put
$$
S_\alpha(\vectx)
 =
\left( S_{(\alpha_i, \alpha_j)}(\vectx) \right)_{1 \le i, j \le l}.
$$
Given a strict partition $\lambda$ of length $l$, we have
\begin{equation}
\label{eq:Schur3}
Q_\lambda(\vectx)
 = 
\begin{cases}
 \Pf S_\lambda(\vectx) &\text{if $l$ is even,} \\
 \Pf S_{\lambda^0}(\vectx) &\text{if $l$ is odd,}
\end{cases}
\end{equation}
where $\lambda^0 = (\lambda_1, \dots, \lambda_l, 0)$.
\end{enumerate}
\end{theorem}

First we show the following stability of Schur $Q$-functions.

\begin{lemma}
\label{lem:stability}
For a strict partition $\lambda$, we have
$$
Q_\lambda(x_1, \dots, x_n, 0)
 =
Q_\lambda(x_1, \dots, x_n).
$$
\end{lemma}

\begin{demo}{Proof}
Let $\vectx = (x_1, \dots, x_n)$ and $l = l(\lambda)$.
Note that $D(x_1, \dots, x_n, 0) = (-1)^n D(x_1, \dots, x_n)$.

If $n+l$ is even, then by definition (\ref{eq:NimmoQ}) we have
$$
Q_\lambda(x_1, \dots, x_n, 0)
 =
\frac{1}{(-1)^n D(\vectx)}
\Pf \begin{pmatrix}
 A(\vectx) & - \vectone_{n,1} & W_\lambda(\vectx) & \vectone_{n,1} \\
 \vectone_{1,n} & 0 & O_{1,l} & 1 \\
 -\trans W_\lambda(\vectx) & O_{l,1} & O_{l,l} & O_{l,1} \\
 -\vectone_{1,n} & -1 & O_{1,l} & 0 
\end{pmatrix},
$$
where $\vectone_{p,q}$ is the all-one matrix of size $p \times q$.
By adding the $(n+1)$st column/row to the last column/row and then expanding the resulting Pfaffian 
along the last column/row, we see that
$$
Q_\lambda(\vectx,0)
 =
\frac{1}{(-1)^n D(\vectx)}
\cdot
(-1)^n \Pf \begin{pmatrix}
 A(\vectx) & W_\lambda(\vectx) \\
 -\trans W_\lambda(\vectx) & O
\end{pmatrix}
 =
Q_\lambda(\vectx).
$$
If $n+l$ is odd, then we have
$$
Q_\lambda(x_1, \dots, x_n, 0)
 =
\frac{1}{ (-1)^n D(\vectx) }
\Pf \begin{pmatrix}
 A(\vectx) & -\vectone_{n,1} & W_\lambda(\vectx) \\
 \vectone_{1,n} & 0 & O_{1,l} \\
 -\trans W_\lambda(\vectx) & O_{l,1} & O_{l,l}
\end{pmatrix}.
$$
By pulling out the common factor $-1$ from the $(n+1)$st row/column and 
then moving the $(n+1)$st row/column to the last row/column, 
we see that
$$
Q_\lambda(\vectx,0)
 =
\frac{1}{ (-1)^n D(\vectx) }
\cdot
(-1)^{l+1}
\Pf \begin{pmatrix}
 A(\vectx) & W_\lambda(\vectx) & \vectone_{n,1} \\
 -\trans W_\lambda(\vectx) & O_{l,l} & O_{l,1} \\
 -\vectone_{1,n} & O_{1,l} & 0
\end{pmatrix}
 =
Q_\lambda(\vectx).
$$
\end{demo}

\begin{demo}{Proof of Theorem~\ref{thm:Schur}}
(1)
By the stability (Lemma~\ref{lem:stability}), we may assume that $n$ is odd. 
Then we have
$$
Q_{(r)}(\vectx)
 =
\frac{ 1 }{ D(\vectx) }
\Pf \begin{pmatrix}
 A(\vectx) & W_{(r)}(\vectx) \\
 - \trans W_{(r)}(\vectx) & 0
\end{pmatrix},
\quad
r \ge 0,
$$
where $W_{(r)}(\vectx)$ is the column vector $\left( \chi(r) x_i^r \right)_{1 \le i \le n}$.
By using 
$$
\sum_{r \ge 0} \chi(r) x_i^r z^r
 =
\frac{ 1 + x_i z }{1 - x_i z },
$$
we see that
$$
\sum_{r \ge 0} Q_{(r)}(\vectx) z^r
 =
\frac{ 1 }{ D(\vectx) }
\Pf \begin{pmatrix}
 A(\vectx) & H_z(\vectx) \\
 -\trans H_z(\vectx) & 0
\end{pmatrix},
$$
where $H_z(\vectx)$ is the column vector with $i$th entry $(1+x_iz)/(1-x_iz)$.
The last Pfaffian is evaluated by using Proposition~\ref{prop:Pf-Schur} with variables $(x_1, \dots, x_n, -1/z)$ 
and we have
$$
\Pf \begin{pmatrix}
 A(\vectx) & H_z(\vectx) \\
 -\trans H_z(\vectx) & 0
\end{pmatrix}
 =
D(\vectx)
\cdot
\prod_{i=1}^n
 \frac{ 1 + x_i z }{ 1 - x_i z }.
$$
This complete the proof of (1).

(2)
By the stability (Lemma~\ref{lem:stability}), we may assume that $n$ is even.
Then we have
$$
Q_{(r,s)}(\vectx)
 =
\frac{ 1 }{ D(\vectx) }
\Pf \begin{pmatrix}
 A(\vectx) & W_{(r)}(\vectx) & W_{(s)}(\vectx) \\
 - \trans W_{(r)}(\vectx) & 0 & 0 \\
 - \trans W_{(r)}(\vectx) & 0 & 0
\end{pmatrix},
\quad r, s \ge 0,
$$
and hence obtain
$$
\sum_{r, s \ge 0} Q_{(r,s)}(\vectx) z^r w^s
 =
\frac{ 1 }{ D(\vectx) }
\Pf \begin{pmatrix}
 A(\vectx) & H_z(\vectx) & H_w(\vectx) \\
 -\trans H_z(\vectx) & 0 & 0 \\
 -\trans H_w(\vectx) & 0 & 0
\end{pmatrix}.
$$
Applying Proposition~\ref{prop:Pf-Schur} with variables $(x_1, \dots, x_n, -1/z, -1/w)$, we see that
$$
\Pf \begin{pmatrix}
 A(\vectx) & H_z(\vectx) & H_w(\vectx) \\
 -\trans H_z(\vectx) & 0 & \dfrac{z-w}{z+w} \\
 -\trans H_w(\vectx) & -\dfrac{z-w}{z+w} & 0
\end{pmatrix}
 =
D(\vectx)
\prod_{i=1}^n
 \frac{ 1 + x_i z }{ 1 - x_i z }
 \frac{ 1 + x_i w }{ 1 - x_i w }
\cdot
\frac{ z - w }{ z + w }.
$$
By splitting the last row/column, we have
\begin{multline*}
\Pf \begin{pmatrix}
 A(\vectx) & H_z(\vectx) & H_w(\vectx) \\
 -\trans H_z(\vectx) & 0 & \dfrac{z-w}{z+w} \\
 -\trans H_w(\vectx) & -\dfrac{z-w}{z+w} & 0
\end{pmatrix}
\\
 =
\Pf \begin{pmatrix}
 A(\vectx) & H_z(\vectx) & H_w(\vectx) \\
 -\trans H_z(\vectx) & 0 & 0 \\
 -\trans H_w(\vectx) & 0 & 0
\end{pmatrix}
+
\Pf \begin{pmatrix}
 A(\vectx) & H_z(\vectx) & \vectzero \\
 -\trans H_z(\vectx) & 0 & \dfrac{z-w}{z+w} \\
 -\trans \vectzero & -\dfrac{z-w}{z+w} & 0
\end{pmatrix}.
\end{multline*}
By expanding the last Pfaffian along the last row/column and using (\ref{eq:Pf-Schur}), we have
$$
\Pf \begin{pmatrix}
 A(\vectx) & H_z(\vectx) & \vectzero \\
 -\trans H_z(\vectx) & 0 & \dfrac{z-w}{z+w} \\
 -\trans \vectzero & -\dfrac{z-w}{z+w} & 0
\end{pmatrix}
 =
\frac{z-w}{z+w} \Pf A(\vectx)
 =
\frac{z-w}{z+w} D(\vectx).
$$
Hence we have
$$
\sum_{r,s \ge 0} Q_{(r,s)}(\vectx) z^r w^s
 =
\frac{z-w}{z+w}
\left(
\prod_{i=1}^n
 \frac{ 1 + x_i z }{ 1 - x_i z }
 \frac{ 1 + x_i w }{ 1 - x_i w }
 -
1
\right).
$$

(3)
By the stability (Lemma~\ref{lem:stability}), we may assume that $n$ is even.
We apply the Pfaffian analogue of the Sylvester identity (Proposition~\ref{prop:Pf-Sylvester}) 
to the matrix $X$ given by
$$
X = \begin{cases}
 \begin{pmatrix}
  A(\vectx) & W_\lambda(\vectx) \\
  - \trans W_\lambda(\vectx) & O_{l,l}
 \end{pmatrix}
 &\text{if $l$ is even,} \\
 \begin{pmatrix}
  A(\vectx) & W_{\lambda^0}(\vectx) \\
  - \trans W_{\lambda^0}(\vectx) & O_{l+1,l+1}
 \end{pmatrix}
 &\text{if $l$ is odd.} \\
\end{cases}
$$
Since $\Pf X( [n] \sqcup \{ n+i, n+j \} )/\Pf X([n]) = Q_{(\lambda_i,\lambda_j)}$ for $i<j$, 
Schur's identity (\ref{eq:Schur3}) immediately follows from Proposition~\ref{prop:Pf-Sylvester}.
\end{demo}

\begin{remark}
We can give a direct proof to the Pfaffian identity (\ref{eq:Schur3}) in the case where $n$ is odd, 
by applying Proposition~\ref{prop:Pf-Sylvester} to the matrix given by
$$
X = \begin{cases}
\begin{pmatrix}
 A(\vectx) & \vectone_{n,1} & W_\lambda(\vectx) \\
 -\vectone_{1,n} & 0 & O_{1,l} \\
 -\trans W_\lambda(\vectx) & O_{l,1} & O_{l,l}
\end{pmatrix}
 &\text{if $l$ is even,}
\\
\begin{pmatrix}
 A(\vectx) & \vectone_{n,1} & W_\lambda(\vectx) & O_{n,1} \\
 -\vectone_{1,n} & 0 & \vectzero_{1,l} & 1 \\
 -\trans W_\lambda(\vectx) & O_{l,1} & O_{l,l} & O_{l,1} \\
 O_{1,n} & -1 & O_{1,l} & 0
\end{pmatrix}
 &\text{if $l$ is odd.}
\end{cases}
$$
\end{remark}

%% file: schurQ5.tex
\section{%
Cauchy-type identity for $Q$-functions
}

In this section, we use the Pfaffian analogue of the Cauchy--Binet formula (Theorem~\ref{thm:Pf-CB}) 
to prove the Cauchy-type identity for Schur $Q$-functions, which corresponds 
to the orthogonality of $Q$-functions.
To prove the Cauchy-type identity, Schur \cite[Abschnitt~IX]{Schur} (see also \cite[\S~4B]{Jozefiak1}) 
used a characterization of $Q_{(n)}$, 
and Macdonald \cite[III.8]{Macdonald} appeal to the theory of Hall--Littlewood functions.
Also the proof given by Hoffman--Humphreys \cite[Chapter~7]{HH} is based 
on the definition of $Q$-functions in terms of vertex operators.
Bijective proofs are given by Worley \cite[Theorem~6.1.1]{Worley} and Sagan \cite[Corollary~8.3]{Sagan}.
Here we give a simple linear algebraic proof.

\begin{theorem}
\label{thm:Cauchy}
(Schur \cite[p.~231]{Schur})
For $\vectx = (x_1, \dots, x_n)$ and $\vecty = (y_1, \dots, y_n)$, we have
\begin{equation}
\label{eq:Cauchy}
\sum_\lambda P_\lambda(\vectx) Q_\lambda(\vecty)
 =
\prod_{i,j=1}^n \frac{ 1 + x_i y_j }{ 1 - x_i y_j},
\end{equation}
where $\lambda$ runs over all strict partitions.
\end{theorem}

The following lemma is obvious, so we omit the proof.

\begin{lemma}
\label{lem:bijection}
Let $n$ be a positive integer and denote by $\Nat$ the set of nonnegative integers.
To a strict partition $\lambda$ we associate the subset $I_n(\lambda) \subset \Nat$ 
given by
$$
I_n(\lambda)
 = 
\begin{cases}
\{ \lambda_1, \dots, \lambda_{l(\lambda)} \} &\text{if $n+l(\lambda)$ is even,} \\
\{ \lambda_1, \dots, \lambda_{l(\lambda)}, 0 \} &\text{if $n+l(\lambda)$ is odd.}
\end{cases}
$$
Then the correspondence $\lambda \mapsto I_n(\lambda)$ gives a bijection 
from the set of all strict partitions to the set of all subsets $I$ of $\Nat$ with $\# I \equiv n \bmod 2$.
\end{lemma}

\begin{demo}{Proof of Theorem~\ref{thm:Cauchy}}
Apply the Pfaffian version of Cauchy--Binet formula (\ref{eq:Pf-CB2}) 
to the matrices
$$
A = A(\vectx),
\quad
B = A(\vecty),
\quad
S = \Bigl( x_i^k \Bigr)_{1 \le i \le n, k \ge 0},
\quad
T = \Bigl( \chi(k) y_i^k \Bigr)_{1 \le i \le n, k \ge 0}.
$$
It follows from the definition of $P$- and $Q$-functions (\ref{eq:NimmoP}) and (\ref{eq:NimmoQ}) 
that for a strict partition $\lambda$ we have
\begin{align*}
P_\lambda(\vectx)
 &=
\frac{(-1)^{\binom{\# I_n(\lambda)}{2}}}{D(\vectx)} 
\Pf \begin{pmatrix}
 A(\vectx) & S([n];I_n(\lambda)) \\
 -\trans S([n];I_n(\lambda)) & O
\end{pmatrix},
\\
Q_\lambda(\vecty)
 &=
\frac{(-1)^{\binom{\# I_n(\lambda)}{2}}}{D(\vecty)} 
\Pf \begin{pmatrix}
 A(\vecty) & T([n];I_n(\lambda)) \\
 -\trans T([n];I_n(\lambda)) & O
\end{pmatrix}.
\end{align*}
Hence, by using Lemma~\ref{lem:bijection} and applying (\ref{eq:Pf-CB2}), we have
\begin{align*}
\sum_\lambda P_\lambda(\vectx) Q_\lambda(\vecty)
 &=
\frac{1}{D(\vectx) D(\vecty)}
\sum_I
 \Pf \begin{pmatrix}
  A(\vectx) & S([n];I) \\
  - \trans S([n];I) & O
 \end{pmatrix}
 \Pf \begin{pmatrix}
  A(\vecty) & T([n];I) \\
  - \trans T([n];I) & O
 \end{pmatrix}
\\
 &=
\frac{ (-1)^{\binom{n}{2}} }
     { D(\vectx) D(\vecty) }
 \Pf \begin{pmatrix}
  A(\vectx) & S \trans T \\
  - T \trans S & -A(\vecty)
 \end{pmatrix},
\end{align*}
where $\lambda$ runs over all strict partitions and $I$ runs over all subsets of $\Nat$ 
with $\# I \equiv n \bmod 2$.
Since the $(i,j)$ entry of $S \trans T$ is given by
$$
\sum_{k \ge 0} x_i^k \cdot \chi(k) y_j^k
 =
\frac{ 1 + x_i y_j }{ 1 - x_i y_j },
$$
we can use the evaluation of the Schur Pfaffian (\ref{eq:Pf-Schur}) 
with variables $(x_1, \dots, x_n, -1/y_1, \dots, \allowbreak -1/y_n)$ 
to obtain
$$
 \Pf \begin{pmatrix}
  A(\vectx) & S \trans T \\
  - T \trans S & -A(\vecty)
 \end{pmatrix}
 =
D(\vectx)
 \cdot
\prod_{i,j=1}^n 
 \frac{ 1 + x_i y_j }{ 1 - x_i y_j}
 \cdot
(-1)^{ \binom{n}{2} } D(\vecty).
$$
This completes the proof.
\end{demo}

%% file: schurQ6.tex
\section{%
Pragacz--J\'ozefiak--Nimmo identity for skew $Q$-functions
}

In this section, we use the Pfaffian analogue of the Cauchy--Binet formula (Theorem~\ref{thm:Pf-CB})
to prove the Pragacz--J\'ozefiak--Nimmo identity for skew $Q$-functions.
Pragacz--J\'ozefiak \cite{JP} and Nimmo \cite{Nimmo} used 
differential operators to prove this Pfaffian identity 
and Stembridge \cite[Theorem~6.2]{Stembridge} gave a combinatorial proof based on the lattice path method.
In the course of our proof, we find a Pfaffian identity which interpolate Nimmo's identity (\ref{eq:NimmoQ}) 
and Schur's identity (\ref{eq:Schur3}).

Skew $Q$-functions $Q_{\lambda/\mu}(x_1, \dots, x_n)$ are uniquely determined by the equation
$$
Q_\lambda(x_1, \dots, x_n, y_1, \dots, y_k)
 =
\sum_\mu Q_{\lambda/\mu}(x_1, \dots, x_n) Q_\mu(y_1, \dots, y_k),
$$
where $\lambda$ is a strict partition and the summation is taken over all strict partitions $\mu$.

\begin{theorem}
\label{thm:JP}
(Pragacz--J\'ozefiak \cite[Theorem~1]{JP}, Nimmo \cite[(2.22)]{Nimmo})
For two sequences $\alpha = (\alpha_1, \dots, \alpha_l)$ and $\beta = (\beta_1, \dots, \beta_m)$ 
of nonnegative integers, let $M_{\alpha/\beta}(\vectx)$ be the $l \times m$ matrix given by
$$
M_{\alpha/\beta}(\vectx)
 =
\Bigl( Q_{(\alpha_i - \beta_{m+1-j})}(\vectx) \Bigr)_{1 \le i \le l, 1 \le j \le m},
$$
where $Q_{(k)}(\vectx) = 0$ for $k < 0$.
For two strict partitions $\lambda$ and $\mu$, we have
\begin{equation}
\label{eq:JP}
Q_{\lambda/\mu}(\vectx)
 =
\begin{cases}
 \Pf \begin{pmatrix} S_\lambda(\vectx) & M_{\lambda/\mu}(\vectx) \\ -\trans M_{\lambda/\mu}(\vectx) & O \end{pmatrix}
 &\text{if $l(\lambda) \equiv l(\mu) \bmod 2$,}
\\
 \Pf \begin{pmatrix} S_\lambda(\vectx) & M_{\lambda/\mu^0}(\vectx) \\ -\trans M_{\lambda/\mu^0}(\vectx) & O \end{pmatrix}
 &\text{if $l(\lambda) \not\equiv l(\mu) \bmod 2$, }
\end{cases}
\end{equation}
\end{theorem}

Note that
$$
\begin{pmatrix} S_\lambda & M_{\lambda/\mu^0} \\ -\trans M_{\lambda/\mu^0} & O \end{pmatrix}
 =
\begin{pmatrix} S_{\lambda^0} & M_{\lambda^0/\mu} \\ -\trans M_{\lambda^0/\mu} & O \end{pmatrix}.
$$

\begin{demo}{Proof}
We denote by $\tilde{Q}_{\lambda/\mu}(\vectx)$ the right hand side of (\ref{eq:JP}) and 
prove
$$
Q_\lambda(\vectx,\vecty)
 =
\sum_\mu \tilde{Q}_{\lambda/\mu}(\vectx) Q_\mu(\vecty).
$$
By the stability (Lemma~\ref{lem:stability}), 
we may assume that the length $l = l(\lambda)$ and the number $k$ of variables in $\vecty$ have 
the same parity.

We apply the Pfaffian analogue of the Cauchy--Binet formula (\ref{eq:Pf-CB1}) to the matrices 
\begin{gather*}
A = S_\lambda(\vectx),
\quad
S = \Bigl( Q_{\lambda_i-r}(\vectx) \Bigr)_{1 \le i \le l, r \ge 0},
\\
B = A(\vecty),
\quad
T = \Bigl( \chi(r) y_i^r \Bigr)_{1 \le i \le k, r \ge 0}.
\end{gather*}
Then, for a strict partition $\mu$, we have
$$
S([l];I_k(\mu))
 =
\begin{cases}
 M_{\lambda/\mu}(\vectx) &\text{if $l(\mu) \equiv k \bmod 2$,} \\
 M_{\lambda/\mu^0}(\vectx) &\text{if $l(\mu) \not\equiv k \bmod 2$,}
\end{cases}
$$
Hence we have
$$
\Pf \begin{pmatrix}
 A(\vectx) & S([l];I_k(\mu)) \\
 -\trans S([l];I_k(\mu)) & O
\end{pmatrix}
 =
\tilde{Q}_{\lambda/\mu}(\vectx).
$$
And it follows from the definition (\ref{eq:NimmoQ}) that
$$
\frac{1}{D(\vecty)}
\Pf \begin{pmatrix}
 A(\vecty) & T([k];I_k(\mu)) \\
 -\trans T([k];I_k(\mu)) & O
\end{pmatrix}
 =
(-1)^{\binom{\# I_k(\mu)}{2}}
Q_\mu(\vecty).
$$
By applying (\ref{eq:Pf-CB1}), we see that
$$
\sum_\mu \tilde{Q}_{\lambda/\mu}(\vectx) Q_\mu(\vecty)
 =
\frac{1}{D(\vecty)}
\Pf \begin{pmatrix}
 S_\lambda(\vectx) & S \trans T \\
 - T \trans S & A(\vecty)
\end{pmatrix}.
$$
Also it follows from the generating function (\ref{eq:Schur1}) of $Q_{(r)}$'s that 
the $(i,j)$ entry of $S \trans T$ is given by
$$
\sum_{r \ge 0} Q_{\lambda_i-r}(\vectx) \cdot \chi(r) y_j^r = Q_{\lambda_i}(\vectx,y_j).
$$
Now we can complete the proof by using the following Theorem.
\end{demo}

\begin{theorem}
\label{thm:NS}
Let $\vectx = (x_1, \dots, x_n)$ and $\vecty = (y_1, \dots, y_k)$ be two sequence of variables.
For a sequence $\alpha = (\alpha_1, \dots, \alpha_l)$ of length $l$, 
let $N_\alpha(\vectx|\vecty)$ be the $l \times k$ matrix defined by
$$
N_\alpha(\vectx|\vecty)
 =
\Bigl(
 Q_{\alpha_i}(\vectx,y_j)
\Bigr)_{1 \le i \le l, 1 \le j \le k}
$$
For a strict partition $\lambda$ of length $l$, we have
\begin{equation}
\label{eq:NS}
Q_\lambda(\vectx,\vecty)
 =
\begin{cases}
\dfrac{ 1 }{ D(\vecty) }
 \Pf \begin{pmatrix}
  S_\lambda(\vectx) & N_\lambda(\vectx|\vecty) \\
  -\trans N_\lambda(\vectx|\vecty) & A(\vecty)
 \end{pmatrix}
 &\text{if $l+k$ is even,} \\
\dfrac{ 1 }{ D(\vecty) }
 \Pf \begin{pmatrix}
  S_{\lambda^0}(\vectx) & N_{\lambda^0}(\vectx|\vecty) \\
  -\trans N_{\lambda^0}(\vectx|\vecty) & A(\vecty)
 \end{pmatrix}
 &\text{if $l+k$ is odd,}
\end{cases}
\end{equation}
where $\lambda^0 = (\lambda_1, \dots, \lambda_l, 0)$.
\end{theorem}

Note that the identity (\ref{eq:NS}) reduces to Nimmo's identity (\ref{eq:NimmoQ}) if $n=0$ 
and to Schur's identity (\ref{eq:Schur3}) if $k=0$.

\begin{demo}{Proof}
We denote by $Q'_\lambda(\vectx|\vecty)$ the right hand side of (\ref{eq:NS}).

First we show that $Q'_\lambda(\vectx|\vecty)$ is stable with respect to $\vecty$, 
that is, 
\begin{equation}
\label{eq:stability'}
Q'_\lambda(\vectx|y_1, \dots, y_k, 0) = Q'_\lambda(\vectx|y_1, \dots, y_k).
\end{equation}
Let $\vecty = (y_1, \dots, y_k)$.
If $l+k$ is even, then we have by using the stability (Lemma~\ref{lem:stability}) 
of $Q_\lambda(\vectx)$, 
$$
Q'_\lambda(\vectx|\vecty,0)
 =
\frac{ 1 }
     { (-1)^k D(\vecty) }
\Pf \begin{pmatrix}
 S_\lambda(\vectx) & T_\lambda(\vectx) & N_\lambda(\vectx|\vecty) & T_\lambda(\vectx) \\
 -\trans T_\lambda(\vectx) & 0 & \vectone_{1,k} & 1 \\
 -\trans N_\lambda(\vectx|\vecty) & -\vectone_{k,1} & A(\vecty) & -\vectone_{k,1} \\
 -\trans T_\lambda(\vectx) & -1 & \vectone_{1,k} & 0
\end{pmatrix},
$$
where $T_\lambda(\vectx)$ is the column vector $\bigl( Q_{\lambda_i}(\vectx) \bigr)_{1 \le i \le l}$.
By adding the $(l+1)$st row/column multiplied by $-1$ to the last row/column 
and then expanding the resulting Pfaffian along the last row/column, we see that
$$
Q'_\lambda(\vectx|\vecty,0)
 =
\frac{ 1 }
     { (-1)^k D(\vecty) }
(-1)^l
\Pf \begin{pmatrix}
 S_\lambda(\vectx) & N_\lambda(\vectx|\vecty) \\
 -\trans N_\lambda(\vectx|\vecty) & A(\vecty)
\end{pmatrix}
 =
Q'_\lambda(\vectx|\vecty).
$$
If $l+k$ is odd, then by moving the last row/column to the $(l+1)$st row/column we have
\begin{align*}
Q'_\lambda(\vectx|\vecty,0)
 &=
\frac{ 1 }
     { (-1)^k D(\vecty) }
\Pf \begin{pmatrix}
 S_\lambda(\vectx) & N_\lambda(\vectx|\vecty) & T_\lambda(\vectx) \\
 -\trans N_\lambda(\vectx|\vecty) & A(\vecty) & -\vectone_{k,1} \\
 -\trans T_\lambda(\vectx) & \vectone_{1,k} & 0
\end{pmatrix}
\\
 &=
\frac{ 1 }
     { (-1)^k D(\vecty) }
(-1)^k
\Pf \begin{pmatrix}
 S_{\lambda^0}(\vectx) & N_{\lambda^0}(\vectx|\vecty) \\
 -\trans N_{\lambda^0}(\vectx|\vecty) & A(\vecty)
\end{pmatrix}
 =
Q'_\lambda(\vectx|\vecty).
\end{align*}

Next we use the Sylvester formula for Pfaffians (Proposition~\ref{prop:Pf-Sylvester}) to prove
\begin{equation}
\label{eq:Schur'}
Q'_\lambda(\vectx|\vecty)
 = 
\begin{cases}
 \Pf S'_\lambda(\vectx|\vecty) &\text{if $l(\lambda)$ is even,} \\
 \Pf S'_{\lambda^0}(\vectx|\vecty) &\text{if $l(\lambda)$ is odd,}
\end{cases}
\end{equation}
where $\lambda^0 = (\lambda_1, \dots, \lambda_{l(\lambda)}, 0)$ and 
the matrix $S'_\alpha(\vectx|\vecty)$ is defined by 
$$
S'_\alpha(\vectx|\vecty)
 =
\Bigl( Q'_{(\alpha_i, \alpha_j)}(\vectx|\vecty) \Bigr)_{1 \le i, j \le l}.
$$
By the stability (\ref{eq:stability'}), we may assume $k$ is even.
In this case the identity (\ref{eq:Schur'}) can be obtained by applying (\ref{eq:Pf-Sylvester2}) 
to the matrix
$$
X
 =
\begin{cases}
\begin{pmatrix}
  S_\lambda(\vectx) & N_\lambda(\vectx|\vecty) \\
  -\trans N_\lambda(\vectx|\vecty) & A(\vecty)
\end{pmatrix}
 &\text{if $l(\lambda)$ is even,} \\
\begin{pmatrix}
  S_{\lambda^0}(\vectx) & N_{\lambda^0}(\vectx|\vecty) \\
  -\trans N_{\lambda^0}(\vectx|\vecty) & A(\vecty)
\end{pmatrix}
 &\text{if $l(\lambda)$ is odd.}
\end{cases}
$$

By comparing two Pfaffian identities (\ref{eq:Schur3}) and (\ref{eq:Schur'}), 
the proof of Theorem~\ref{thm:NS} is reduced to showing 
\begin{gather}
\label{eq:Q=Q'1}
Q_{(r)}(\vectx,\vecty) = Q'_{(r)}(\vectx|\vecty),
\\
\label{eq:Q=Q'2}
Q_{(r,s)}(\vectx,\vecty) = Q'_{(r,s)}(\vectx|\vecty).
\end{gather}
We prove these equality by considering the generating functions.
If we put
$$
F_z(u_1, \dots, u_m) = \prod_{i=1}^m \frac{ 1 + u_i z }{ 1 - u_i z},
$$
then by virtue of (\ref{eq:Schur1}) and (\ref{eq:Schur2}) 
the identities (\ref{eq:Q=Q'1}) and (\ref{eq:Q=Q'2}) follow from
\begin{equation}
\label{eq:F=F'1}
\sum_{r \ge 0} Q'_{(r)}(\vectx|\vecty) z^r
 =
F_z(\vectx,\vecty),
\end{equation}
and
\begin{equation}
\label{eq:F=F'2}
\sum_{r,s \ge 0} Q'_{(r,s)}(\vectx|\vecty) z^r w^s
 =
\frac{z-w}{z+w} \left( F_z(\vectx,\vecty) F_w(\vectx,\vecty) - 1 \right),
\end{equation}
respectively.

By the stability (\ref{eq:stability'}) we may assume $k$ is odd for the proof of (\ref{eq:F=F'1}).
If $k$ is odd, then
$$
Q'_{(r)}(\vectx|\vecty)
 =
\frac{ 1 }{ D(\vecty) }
\Pf \begin{pmatrix}
 0 & N_{(r)}(\vectx|\vecty) \\
 -\trans N_{(r)}(\vectx|\vecty) & A(\vecty)
\end{pmatrix}
$$
Since $\sum_{r \ge 0} Q(\vectx,y_j) z^j = F_z(\vectx) (1 + y_j z)/(1 - y_j z)$ by (\ref{eq:Schur1}), 
we have
\begin{align*}
\sum_{r \ge 0} Q'_{(r)}(\vectx|\vecty) z^r
 &=
\frac{ 1 }{ D(\vecty) }
\Pf \begin{pmatrix}
 0 & F_z(\vectx) \trans H_z(\vecty) \\
 - F_z(\vectx) H_z(\vectx) & A(\vecty)
\end{pmatrix}
\\
 &=
\frac{ 1 }{ D(\vecty) }
F_z(\vectx)
\Pf \begin{pmatrix}
 0 & \trans H_z(\vecty) \\
 - H_z(\vecty) & A(\vecty)
\end{pmatrix},
\end{align*}
where $H_z(\vecty)$ is the column vector $\left( (1+y_iz)/(1-y_iz) \right)_{1 \le i \le k}$.
By applying Proposition~\ref{prop:Pf-Schur} with variables $(-1/z, y_1, \dots, y_n)$, 
we have
$$
\Pf \begin{pmatrix}
 0 & \trans H_z(\vecty) \\
 - H_z(\vecty) & A(\vecty)
\end{pmatrix}
 =
D(\vecty) F_z(\vecty),
$$
and obtain (\ref{eq:F=F'1}).

By the stability (\ref{eq:stability'}) we may assume $k$ is even for the proof of (\ref{eq:F=F'2}).
If $k$ is even, then
$$
Q'_{(r,s)}(\vectx|\vecty)
 =
\frac{1}{D(\vecty)}
\Pf \begin{pmatrix}
 0 & Q_{(r,s)}(\vectx) & N_{(r)}(\vectx|\vecty) \\
 -Q_{(r,x)}(\vectx) & 0 & N_{(s)}(\vectx|\vecty) \\
 -\trans N_{(r)}(\vectx|\vecty) & -\trans N_{(s)}(\vectx|\vecty) & A(\vecty)
\end{pmatrix}
$$
for $r$, $s \ge 0$.
Hence it follows from (\ref{eq:Schur1}) and (\ref{eq:Schur2}) that
$$
\sum_{r, s \ge 0} Q'_{(r,s)}(\vectx|\vecty) z^r w^s
 =
\frac{1}{D(\vecty)}
\Pf \begin{pmatrix}
 0 & G_{z,w}(\vectx) & F_z(\vectx) \trans H_z(\vecty) \\
 - G_{z,w}(\vectx) & 0 & F_w(\vectx) \trans H_w(\vecty) \\
 - F_z(\vectx) H_z(\vecty) & - F_w(\vectx) H_w(\vecty) & A(\vecty)
\end{pmatrix},
$$
where
$$
G_{z,w}(\vectx) = \frac{z-w}{z+w} (F_z(\vectx) F_w(\vectx) - 1).
$$
By splitting the first row/column and then pulling out the common factor $F_z(\vectx)$ and $F_w(\vectx)$ 
from the $1$st and $2$nd rows/columns, we see that
\begin{align*}
\sum_{r, s \ge 0} Q'_{(r,s)}(\vectx|\vecty) z^r w^s
 &=
\frac{ F_z(\vectx) F_w(\vectx) }
     { D(\vecty) }
\Pf \begin{pmatrix}
 0 & \dfrac{z-w}{z+w} & \trans H_z(\vecty) \\
 -\dfrac{z-w}{z+w} & 0 & \trans H_w(\vecty) \\
 -\trans H_z(\vecty) & -\trans H_w(\vecty) & A(\vecty)
\end{pmatrix}
\\
 &\quad
-
\frac{ 1 }
     { D(\vecty) }
\Pf \begin{pmatrix}
 0 & \dfrac{z-w}{z+w} & O_{1,k} \\
 - \dfrac{z-w}{z+w} & 0 & \trans H_w(\vecty) \\
 O_{k,1} & -\trans H_w(\vecty) & A(\vecty)
\end{pmatrix}.
\end{align*}
The first Pfaffian is evaluated by using the Schur Pfaffian with $(-1/z,-1/w,y_1, \dots, y_n)$ 
and we see that
$$
\Pf \begin{pmatrix}
 0 & \dfrac{z-w}{z+w} & \trans H_z(\vecty) \\
 -\dfrac{z-w}{z+w} & 0 & \trans H_w(\vecty) \\
 -\trans H_z(\vecty) & -\trans H_w(\vecty) & A(\vecty)
\end{pmatrix}
 =
\frac{z-w}{z+w}
F_z(\vecty) F_w(\vecty)
D(\vecty)
$$
By expanding the second Pfaffian along the first column/row, we see that it equals to
$$
\Pf \begin{pmatrix}
 0 & \dfrac{z-w}{z+w} & O_{1,k} \\
 - \dfrac{z-w}{z+w} & 0 & \trans H_w(\vecty) \\
 O_{k,1} & -\trans H_w(\vecty) & A(\vecty)
\end{pmatrix}
 =
\frac{z-w}{z+w} D(\vecty).
$$
Therefore we obtain
$$
\sum_{r, s \ge 0} Q'_{(r,s)}(\vectx|\vecty) z^r w^s
 =
\frac{z-w}{z+w}
\left( F_z(\vectx,\vecty) F_w(\vectx,\vecty) - 1 \right).
$$
This completes the proof of Theorem~\ref{thm:NS} and hence Theorem~\ref{thm:JP}.
\end{demo}

%% file: schurQ7.tex
\section{%
Littlewood-type identity for $Q$-functions
}

In this section, we prove the following Littlewood-type identity for $Q$-functions. 
This identity is a special case ($t=\sqrt{-1}$) of \cite[(1.21)]{Kawanaka2} for Hall--Littlewood functions, 
which is essentially proved in \cite{Kawanaka1} by using the representation theory 
of finite Chevalley groups.

\begin{theorem}
\label{thm:Littlewood}
(Kawanaka \cite{Kawanaka2})
For $\vectx = (x_1, \dots, x_n)$, we have
\begin{equation}
\label{eq:Littlewood}
\sum_\lambda
 \bigl( 1+\sqrt{-1} \bigr)^{l(\lambda)}
 P_\lambda(\vectx)
 =
\prod_{i=1}^n \frac{1 + \sqrt{-1} x_i}{1 - x_i}
\prod_{1 \le i < j \le n} \frac{1 + x_i x_j}{1 - x_i x_j},
\end{equation}
where $\lambda$ runs over all strict partitions.
\end{theorem}

\begin{remark}
The right hand side of (\ref{eq:Littlewood}) is one of the simplest example of products 
involving the factor $\prod_{1 \le i < j \le n} (1+x_ix_j)/(1-x_ix_j)$ 
that is a(n infinite) linear combination of $P$- or $Q$-functions.
Recall that a symmetric polynomial $f(\vectx)$ is a linear combination 
of Schur $Q$-functions if and only if $f(t, -t, x_3, \dots, x_n)$ is independent of $t$.
(See \cite[III (8.5)]{Macdonald} for example.)
Consider a symmetric power series of the form
$$
f_n(\vectx) =
\prod_{i=1}^n
 \frac{ \prod_{j=1}^r (1 - \alpha_j x_i) }
      { \prod_{j=1}^s (1 - \beta_j x_i) }
\prod_{1 \le i < j \le n}
 \frac{ 1 + x_i x_j }{ 1 - x_i x_j },
$$
where $\{ \alpha_1, \dots, \alpha_r \} \cap \{ \beta_1, \dots, \beta_s \} = \emptyset$.
Then 
$$
f_n(t,-t, x_3, \dots, x_n)
 =
\frac{ \prod_{j=1}^r (1 - \alpha_j t) (1 + \alpha_j t) }
     { \prod_{j=1}^s (1 - \beta_j t) (1 + \beta_j t) }
\cdot
\frac{ 1 - t^2 }
     { 1 + t^2 }
\cdot
f_{n-2}(x_3, \dots, x_n)
$$
is independent of $t$ if and only if 
$$
\{ \alpha_1, \dots, \alpha_r, -\alpha_1, \dots, -\alpha_r, 1, -1 \}
 =
\{ \beta_1, \dots, \beta_s, -\beta_1, \dots, \beta_s, \sqrt{-1}, -\sqrt{-1} \}
$$
as multisets.
Thus $f_n(\vectx)$ is an infinite linear combination of $P$-functions if and only if 
$r = s$ and $\alpha_1 = \pm \sqrt{-1}$, $\beta_1 = \pm 1$ 
and $\alpha_k = - \beta_k$ for $2 \le k \le r$
up to permutation of $\alpha_1, \dots, \alpha_r$ and $\beta_1, \dots, \beta_r$.
\end{remark}

\begin{demo}{Proof of Theorem~\ref{thm:Littlewood}}
By the stability (Lemma~\ref{lem:stability}), we may assume $n$ is even.
We apply the Pfaffian version of the minor-summation formula (Theorem~\ref{thm:Pf-IW}) to
the matrices
$$
A = A(\vectx),
\quad
S = \Bigl( x_i^k \Bigr)_{1 \le i \le n, k \ge 0},
$$
and the skew-symmetric matrix $B$ whose $(i,j)$ entry, $0 \le i < j$, is given by
$$
B_{ij} =
\begin{cases}
 -\alpha &\text{if $i=0$,} \\
 -\alpha^2 &\text{if $i>0$,}
\end{cases}
$$
where $\alpha = 1+\sqrt{-1}$.

By using (\ref{eq:Pf(-X)}), (\ref{eq:Pf-expansion}) and the induction on $\# I$, 
we see that the subpfaffian $\Pf B(I)$ of $B$ corresponding to a even-element subset $I \subset \Nat$ is given by
$$
\Pf B(I)
 =
\begin{cases}
(-1)^{\binom{\# I}{2}} \alpha^{\# I - 1} &\text{if $0 \in I$,} \\
(-1)^{\binom{\# I}{2}} \alpha^{\# I} &\text{if $0 \not\in I$.}
\end{cases}
$$
Since $n$ is even, strict partitions $\lambda$ are in bijection with even-element subsets $I_n(\lambda)$ 
of $\Nat$ by Lemma~\ref{lem:bijection} and 
$$
\Pf B(I_n(\lambda)) = (-1)^{\binom{\# I_n(\lambda)}{2}} \alpha^{l(\lambda)}.
$$
Also it follows from Nimmo's identity (\ref{eq:NimmoP}) that
$$
\Pf \begin{pmatrix}
 A & S([n];I_n(\lambda)) \\
 -\trans S([n];I_n(\lambda)) & O
\end{pmatrix}
 =
(-1)^{\binom{\# I_n(\lambda)}{2}}
D(\vectx) P_\lambda(\vectx).
$$
Therefore we have
$$
\Pf (A - S B \trans S)
 =
\sum_I
 \Pf B(I) 
 \Pf \begin{pmatrix}
  A & S([n];I) \\
  -\trans S([n];I) & O
 \end{pmatrix}
 =
D(\vectx)
\sum_\lambda \alpha^{l(\lambda)} P_\lambda(\vectx),
$$
where $I$ runs over all even-element subsets of $\Nat$ 
and $\lambda$ runs over all strict partitions.

By direct computations, we see that the $(i,j)$-entry of $S B \trans S$ is equal to
$$
\sum_{k,l \ge 0} b_{k,l} x_i^k x_j^l
 =
- \alpha \frac{x_j-x_i}{(1-x_i)(1-x_j)}
- \alpha^2 \frac{ x_i x_j (x_j - x_i) }{(1-x_i x_j)(1-x_i)(1-x_j)},
$$
and the $(i,j)$ entry of $A - S B \trans S$ is equal to
$$
\frac{x_j-x_i}{x_j+x_i} - \sum_{k,l \ge 0} b_{k,l} x_i^k x_j^l
 =
\frac{1 + \sqrt{-1} x_i}{1 - x_i}
\frac{1 + \sqrt{-1} x_j}{1 - x_j}
\frac{(1+x_ix_j)(x_j-x_i)}{(1-x_ix_j)(x_j+x_i)}.
$$
Hence, by using Proposition~\ref{prop:Pf-Schur} with variables $(x_1-1/x_1, \dots, x_n-1/x_n)$, 
we have
\begin{align*}
\Pf ( A - S B \trans S )
 &=
\prod_{i=1}^n
\frac{1 + \sqrt{-1} x_i}{1 - x_i}
\Pf \left(
 \frac{ (x_j-x_i)(1+x_ix_j) }
      { (x_j+x_i)(1-x_ix_j) }
\right)_{1 \le i, j \le n}
\\
 &=
\prod_{i=1}^n
\frac{1 + \sqrt{-1} x_i}{1 - x_i}
\prod_{1 \le i < j \le n}
 \frac{ (x_j-x_i)(1+x_ix_j) }
      { (x_j+x_i)(1-x_ix_j) }.
\end{align*}
This completes the proof.
\end{demo}

By considering the real and imaginary parts of Theorem~\ref{thm:Littlewood}, we obtain 

\begin{corollary}
If we put
$$
a_l = \begin{cases}
 (-1)^k 2^{2k}       &\text{if $l = 4k$,} \\
 (-1)^k 2^{2k}   &\text{if $l = 4k+1$,} \\
 0                   &\text{if $l = 4k+2$,} \\
 (-1)^{k+1} 2^{2k+1} &\text{if $l = 4k+3$,} \\
\end{cases}
\quad
b_l = \begin{cases}
 0                   &\text{if $l = 4k$,} \\
 (-1)^k 2^{2k}   &\text{if $l = 4k+1$,} \\
 (-1)^k 2^{2k+1}     &\text{if $l = 4k+2$,} \\
 (-1)^k 2^{2k+1} &\text{if $l = 4k+3$,} \\
\end{cases}
$$
then we have
\begin{gather*}
\sum_\lambda a_{l(\lambda)} P_\lambda(\vectx)
 =
\frac{ 1 - e_2 + e_4 - e_6 + \dots }
     { 1 - e_1 + e_2 - e_3 + \dots }
\prod_{1 \le i < j \le n}
 \frac{1 + x_i x_j}{1 - x_i x_j},
\\
\sum_\lambda b_{l(\lambda)} P_\lambda(\vectx)
 =
\frac{ e_1 - e_3 + e_5 - \dots }
     { 1 - e_1 + e_2 - e_3 + \dots }
\prod_{1 \le i < j \le n}
 \frac{1 + x_i x_j}{1 - x_i x_j},
\end{gather*}
where $e_k = e_k(\vectx)$ is the $k$th elementary symmetric polynomial.
\end{corollary}

\begin{demo}{Proof}
Since we have
$$
\bigl( 1 + \sqrt{-1} \bigr)^l
 =
a_l + b_l \sqrt{-1},
$$
we obtain this corollary from Theorem~\ref{thm:Littlewood}.
\end{demo}